\documentclass[a4paper,12pt,twoside]{amsart}
\usepackage{amsmath,amsthm,amsfonts}
\usepackage[T1]{fontenc}

\newtheorem{thm}{Theorem}[section]

\newtheorem{lem}[thm]{Lemma}

\newtheorem{prop}[thm]{Proposition}
\newtheorem{proposition}[thm]{Proposition}

\newtheorem{defi}[thm]{Definition}
\newtheorem{nota}[thm]{Notation}

\theoremstyle{definition}

\newtheorem{rem}{Remark}
\newtheorem{rems}{Remarks}

\newcommand{\Cal}{\mathcal}
\newcommand{\cal} {\mathcal}

\newcommand{\R}{{\mathbb{R}}}

\newcommand{\E}{{\mathbb{E}}}

\newcommand{\Z}{{\mathbb{Z}}}
\newcommand{\N}{{\mathbb{N}}}
\newcommand{\T}{{\mathbb{T}}}
\newcommand{\PP}{{\mathbb{P}}}

\def \RR {\R}

\def \EE {\E}

\def \ZZ {\Z}

\newcommand{\Sp}{{\mathbb{S}}}

\def \ev {\cal E}
\def \ad {\cal D}
\def \ead {\cal E\cal D}

\def \Card {{\rm Card}}

\def \mod {{\rm \ mod \ }}
\def \for {{\rm \ for \ }}
\def \ae {{\rm -a.e. \ }}
\def \ifrm {{\rm if \ }}

\def\eop{\qed}
\def\proof {\vskip -3mm {{\it Proof}.}}

\title [Limit directions of a vector cocycle]
{Limit directions of a vector cocycle,\\ remarks and examples}

\bigskip
\date {April 2014}
\author{Jean-Pierre Conze}
\author{Stéphane Le Borgne}
\address{IRMAR, UMR CNRS 6625, Universit\'e de Rennes I
\vskip 0mm Campus de Beaulieu, 35042 Rennes Cedex, France}
\email{conze@univ-rennes1.fr}
\email{stephane.leborgne@univ-rennes1.fr}


\begin{document}
\baselineskip 15pt \parindent=0mm

\maketitle \centerline {IRMAR, UMR CNRS 6625,} \centerline {
University of Rennes I, France}

\begin{abstract} We study the set $\ad(\Phi)$ of limit directions
of a vector cocycle $(\Phi_n)$ over a dynamical system, i.e., the
set of limit values of $\Phi_n(x) /\|\Phi_n(x)\|$ along subsequences
such that $\|\Phi_n(x)\|$ tends to $\infty$. This notion is natural
in geometrical models of dynamical systems where the phase space is
fibred over a basis with fibers isomorphic to $\R^d$, like systems
associated to the billiard in the plane with periodic obstacles. It
has a meaning for transient or recurrent cocycles.

Our aim is to present some results in a general context as well as
for specific models for which the set of limit directions can be
described. In particular we study the related question of sojourn in
cones of the cocycle when the invariance principle is satisfied.
\end{abstract}

\tableofcontents

\section*{Introduction}

Let $(X,\mu, T)$ be an ergodic dynamical system and $\Phi$ be a
measurable function on $X$ with values in $\mathbb{R}^d$. The
ergodic sums $\Phi_n(x) := \sum_{k=0}^{n-1}\Phi(T^kx), n \geq 1$,
define a vector process. When $\Phi$ is integrable and not centered,
this process tends a.s. to $\infty$ in the direction of the mean
$\int \Phi \, d\mu$. A general question, when $\Phi$ is centered or
for a measurable non integrable $\Phi$, is to find in which
directions at infinity the ergodic sums are going. The set of these
directions is a kind of boundary for the cocycle $(\Phi_n)$, i.e.
for the process of ergodic sums.

This leads to the notion of limit directions and to the
cohomologically invariant notion of essential limit directions. The
limit directions of a vector cocycle $(\Phi_n)$ over a dynamical
system can be defined as the limit values of $\Phi_n(x)
/\|\Phi_n(x)\|$ along subsequences such that $\|\Phi_n(x)\|$ tends
to $\infty$.

The notion of limit directions is natural in geometrical models of
dynamical systems where the phase space is fibred over a basis with
fibers isomorphic to $\R^d$, like the dynamical systems associated
to the billiard in the plane with periodic obstacles. It has a
meaning for recurrent cocycles as well as for transient cocycles.

Our aim is to present some results in a general context (Section
\ref{asympDirec}) and for specific models where the set of limit
directions can be made explicit. In Subsection \ref{oscildim1}
1-dimensional cocycles are considered and some classical results are
recalled or slightly extended.

In Section \ref{chaotMap}, we apply properties like the CLT for
subsequences or the invariance principle to study essential limit
directions and the behavior of the process induced by the cocycle on
the sphere. For $d \geq 2$, when $\Phi$ satisfies a Central Limit
Theorem, one can think that the limit behavior of the sums is
analogous to that of a Brownian motion, in particular in terms of
visit of cones. In the last subsection \ref{invprincSect} this is
shown to be the case, at least if $\Phi$ satisfies Donsker's
invariance principle.

\section{\bf Preliminaries \label{prelim}}

Let $(X, {\cal B}, \mu, T)$ be a dynamical system, where $(X, {\cal
B})$ is a standard Borel space, $T$ an invertible measurable map $T:
X \rightarrow X$ and $\mu$ a probability measure which is
$T$-invariant. Let $\Phi$ be a measurable function on $X$ with
values in $G =\RR^d$. The process $(\Phi \circ T^n)_{n \geq 1}$ is
stationary. Recall that, under general assumptions, every stationary
process $(X_n)$ with values in $\R^d$ can be represented in such a
way for some dynamical system and some measurable $\Phi$.

Part of the results below are valid when $\mu$ is only supposed to
be a $\sigma$-finite $T$-quasi-invariant measure such that $T$ is
{\it conservative} for $\mu$, i.e., for every measurable $B$ in $X$,
for $\mu$-a.e. $x \in B$, there is $n(x) > 0$ such that $T^{n(x)} x
\in B$. Nevertheless for the sake of simplicity, excepted in
Subsection \ref{oscildim1}, we will restrict the presentation to the
framework of a probability invariant measure $\mu$.

Excepted in Section \ref{oscildim1}, the system $(X, \mu, T)$ is
supposed to be ergodic. All equalities are understood to hold
$\mu$-a.e. All sets that we consider are measurable and (unless the
contrary is explicitly stated) with positive measure.

To $\Phi$ is associated a {\it cocycle} $(\Phi_n)_{n \in \N}$
defined by $\Phi_0 (x) = 0$,
$$\Phi_n (x) = \Phi (x) + \, ... \, + \Phi (T^ {n-1} x), \text{ for } n \ge 1,
\text{ and }  \Phi_n = -\Phi_{-n} \circ T^{n}, \text{ for } n < 0,$$
and a map $T_\Phi$ (called {\it skew product}) acting on $X\times
\RR^d$ by
\begin{equation}
T_\Phi: (x, y)\rightarrow (T x,y + \Phi (x)).\label {cylind}
\end{equation}

The cocycle relation $\Phi_{n+p}(x) = \Phi_n(x) + \Phi_p(T^n x),
\forall n,p \in \Z$ is satisfied. The cocycle gives the position in
the fiber after $n$ iterations of $T_\Phi$:
$$T_\Phi^n(x,y) = (T^n x, \, y + \Phi_n(x)).$$
The cocycle $(\Phi_n)$ can be viewed as a {\it "stationary" walk} in
$\RR^d$ "driven" by the dynamical system $(X, \mu ,T)$. It is also
the sequence of ergodic sums of $\Phi$ for the action of $T$. We
will use as well the notation $(\Phi, T)$.

The Lebesgue measure on $\RR^d$ is denoted by $m(dy)$ or simply
$dy$. The map $T_\Phi$ leaves invariant the measure $\mu \times m$
denoted by $\lambda$.

Recall that a cocycle $(\Phi_n)$ over $(X, \mu, T)$ is {\it
transient} if $\lim_n \|\Phi_n(x)\| = +\infty$, for a.e. $x \in X$.
It is recurrent if $\liminf_n \|\Phi_n(x)\| < \infty$, for a.e. $x
\in X$. It is well known that, when $T$ is ergodic, a cocycle is
either transient or recurrent (see the comment below).

Recurrence of the cocycle is equivalent to conservativity of the map
$T_\Phi$ for the measure $\lambda$. When $(\Phi_n)$ is recurrent,
then $(\Phi_n(x))$ returns for a.e. $x$ infinitely often in any
neighborhood of the origin. In dimension 1, if $\Phi$ is integrable
and $(X, \mu, T)$ is ergodic, then $(\Phi_n)_{n\in\ZZ}$ is recurrent
if and only if $\mu(\Phi) = 0$. In higher dimension, recurrence
requires stronger assumptions.

\vskip 3mm {\bf Induced map}

Let us recall some definitions and notations about induced maps.

Let $B$ be a measurable set of positive $\mu$-measure. On $B$
equipped with the measure $\mu_B = \mu(B)^{-1} \mu_{|B}$, the
induced transformation is $T_B(x) = T^{R(x)}(x)$, where $R(x)$ is
the return time $R(x) = R_B(x) := \inf\{j \ge 1: T^j x \in B\}$. The
return time is well defined for a.e. $x \in B$ by conservativity of
the system. We induce\footnote{In short the function $\Phi$ itself
will also be called "cocycle" and $\Phi^B$ "induced cocycle" on
$B$.} $\Phi$ on $B$ by putting
\begin{eqnarray}
\Phi^B(x) :=\Phi_{R(x)}(x)= \sum_{j= 0}^{R(x)- 1} \Phi(T^j
x).\label{inducCoc}
\end{eqnarray} The "induced" cocycle is $\Phi_n^B(x):=
\Phi^B(x) + \Phi^B(T_B x)\,\cdots \, + \Phi^B(T_B^{n-1}x)$, for $n
\ge 1$.

If $\Phi$ is recurrent, then each induced cocycle $(\Phi_n^B)$ is
recurrent. Indeed $(T_B)_{\Phi_B}$ is the induced map on $B\times G$
of $T_\Phi$ which is conservative.

If $T$ is ergodic, then $(B, \mu_B, T_B)$ is ergodic. The converse
is true when $X = \bigcup_n T^n B$.

When the map $T$ is ergodic, the above formulas can be extended to
$X$ by setting, for every measurable set $B$ of positive measure,
for a.e. $x \in X$:
\begin{eqnarray}
R_{B}(x) = \inf\{j \geq 1: T^j x \in B \}, \ \  \Phi^{B}(x) =
\sum_{j=0}^{R_B(x) -1} \Phi(T^j x). \label{passage}
\end{eqnarray}

Recall that two $\RR^d$-valued cocycles $(\Phi^1, T)$ and
$(\Phi^2,T)$ over the dynamical system $(X, \mu, T)$ are $\mu$-{\it
cohomologous} with {\it transfer function} $\Psi$, if there is a
measurable map $\Psi : X \rightarrow \RR^d$ such that
\begin{equation}
\Phi^1(x)  = \Phi^2(x) + \Psi(T x) - \Psi(x), \text{ a.e. }
\label{cohom}
\end{equation}
$\Phi$ is a {\it $\mu$-coboundary}, if it is cohomologous to 0.

We choose a norm $\| \ \|$ on $\R^d$. We will use the inequality
\begin{eqnarray}
|\|\Phi_{n+1}(x)\| -  \|\Phi_n(T x)\|| \leq \|\Phi(x)\|.
\label{inegTriang}
\end{eqnarray}

\vskip 3mm
\section{\bf Limit directions of a vector cocycle, general properties \label{asympDirec}}

\vskip 3mm {\bf "0\,-1" properties for a cocycle}

Let $(\Phi_n)$ be a cocycle over an ergodic dynamical system $(X,
\mu, T)$. Some of its limit properties are related to the ergodicity
of the skew product $T_\Phi$. For example, equirepartition
properties (comparison of the number of visits to sets of finite
measure) are given by the ratio ergodic theorem when the skew
product $T_\Phi$ is ergodic.

There are also limit properties which do not a priori require
ergodicity of the skew product, but appear as "0\,-1" properties, in
the sense that either they are satisfied by a.e. $x$, or are not
satisfied by a.e. $x$.

More precisely, let ${\cal P}(x)$ be a property which, for $x \in
X$, is satisfied or not by the sequence $(\Phi_n(x))$. If the set $
{\cal A}_{\cal P} := \{x: {\cal P(x)} {\rm \ is \ true}\}$ is
measurable and invariant by the map $T$, then by ergodicity of $(X,
\mu, T)$ this set has measure 0 or 1: either ${\cal P}(x)$ is true
for a.e. $x$, or ${\cal P}(x)$ is false for a.e. $x$.

Sometimes, for an asymptotic property ${\cal P}$, the set ${\cal
A}_{\cal P}$ can be described in term of $\limsup$ of a sequence of
sets and its invariance by the map $T$ can easily be checked.

The dichotomy between {\it recurrence and transience} of a cocycle
is an example of a "0\,-1" property: the property ${\cal R}$ "the
cocycle is recurrent" corresponds to the set ${\cal A}_{\cal R} =
\bigcup_{M \geq 1} \bigcap_{N \geq 1} \bigcup_{n \geq N} A_n^M$,
where $A_n^M = \{x: \|\Phi_n(x)\| \leq M \}$.

Indeed, from the inequality (\ref{inegTriang}) it follows $T^{-1}
{\cal A}_{\cal R} = {\cal A}_{\cal R}$. Therefore, when $(X, \mu,
T)$ is ergodic, either for $\mu$ a.e. every $x$, $\lim_n
\|\Phi_n(x)\| = +\infty$, or for a.e. $x$ the cocycle $(\Phi_n(x))$
returns infinitely often in some compact set depending on $x$. In
the latter case, an argument based on Poincaré recurrence lemma
implies that the cocycle returns to any neighborhood of 0, for a.e.
$x$.

We give below another example: the notion of limit direction.

\subsection{\bf Limit directions}

\

\vskip 3mm {\bf Essential values and regularity}

First we recall the classical notion of essential values of a
recurrent cocycle with values in an abelian lcsc group $G$ (cf. K.
Schmidt \cite{Sc77}). A point $\infty$ is added to $G$ with the
natural notion of neighborhood. For our purpose, we restrict
ourselves to the case $G = \R^d$.

\begin{defi}\label{valess} {\rm An element $a \in G \cup \{\infty\}$
is {\it an essential value} of the cocycle $(\Phi, T)$ (with respect
to $\mu$) if, for every neighborhood $V(a)$ of $a$, for every
measurable subset $B$ of positive measure,
\begin{eqnarray}
\mu(B\cap T^{-n} B \cap \{x: \Phi_n(x) \in V(a) \}\bigr) > 0, {\rm \
for \ some \ } n \geq 0. \label{visitVa}
\end{eqnarray} } \end{defi}

The property (\ref{visitVa}) can be stated in the equivalent way:
\begin{eqnarray}
\mu (\{x \in B: \Phi_n^B(x) \in V(a)\}) > 0, {\rm \ for \ some \ } n
\geq 0. \label{visitVa1}
\end{eqnarray}

We denote by ${\overline \ev}(\Phi)$ the set of essential values of
the cocycle $(\Phi, T)$ and by $\ev(\Phi) = {\overline
\ev}(\Phi)\cap G$ the {\it set of finite essential values}.

Let us recall some facts. The set $\ev(\Phi)$ is a closed subgroup
of $G$. A cocycle $\Phi$ is a coboundary if and only if ${\overline
{\cal E}}(\Phi) =\{ 0 \}$. We have ${\cal E}(\Phi {\mod {\cal
E}}(\Phi)) = \{0\}$. Two cohomologous cocycles have the same set of
essential values.

It is well known (\cite{Sc77}, \cite{Aa97}) that the set $\ev(\Phi)$
coincides with ${\mathcal P}(\Phi)$, the group of periods $p$ of the
measurable $T_\Phi$-invariant functions on $X \times G$, i.e., the
elements $p \in G$ such that for every $T_\Phi$-invariant $F$,
$F(x,y +p) = F(x,y), \lambda- a.e.$ This shows that ${\cal E}(\Phi)
= G$ if and only if $(X \times G, \lambda_\chi, T_\Phi)$ is ergodic.

\begin{defi}\label{regul0} {\rm We say that the cocycle defined
by $\Phi$ is {\it regular}, if it is cohomologous to a cocycle which
has values in a closed subgroup $H$ of $G$ and is ergodic on $X
\times H$. The group $H$ in the definition is ${\cal
E}(\Phi)$.}\end{defi}

Now we consider the notion of limit directions and essential limit
directions. The cocycle can be recurrent or transient.

\vskip 3mm
\goodbreak {\bf Limit directions}

For $v \in \R^d \setminus \{0\}$, let $\tilde v := v/\|v\|$ be the
corresponding unit vector in the unit sphere $\Sp^{d-1}$. For every
$\RR^d$-valued cocycle $(\Phi_n)$, we obtain a process (directional
process) $(\tilde \Phi_n)_{n \geq 1}$ with values in $\Sp^{d-1}$
(defined outside the values $(n, x)$ such that $\Phi_n(x) = 0$).

\begin{defi} \label{asymp_dir_def} {\rm A vector $u$ is {\it a limit direction}
of the cocycle $(\Phi_n(x))$ at $x$, if there exists a subsequence
$(n_k(x))$ such that $\|\Phi_{n_k}(x)\| \to \infty$ and
$\Phi_{n_k}(x)/ \|\Phi_{n_k}(x)\|$ converges to $u$.

The subset for which the property ${\cal P}_u$: "$u$ is a limit
direction of $(\Phi_n(x))$" holds is
\begin{eqnarray}
{\cal A}(u) = \bigcap_{V,M} \bigcap_N \bigcup_{n \geq N} \{x \in X:
\|\Phi_n(x)\| > M  {\rm \ and \ } \Phi_n(x)/\|\Phi_n(x)\| \in V\}.
\label{Asymp-u}
\end{eqnarray}
where the intersection is taken over a countable basis of
neighborhoods $V$ of $u$ and the positive integers $M$.

The set of limit directions of the cocycle $(\Phi_n(x))$ for $x \in
X$ is defined as
$$\ad(\Phi)(x) := \{u: \exists \  (n_k(x)): \|\Phi_{n_k(x)}(x)\|
\to \infty {\rm \ and \ } \Phi_{n_k(x)}(x)/ \|\Phi_{n_k(x)}(x)\| \to
u\}.$$ }\end{defi}

\begin{rems} \label{remderiv} {\rm a) From (\ref{inegTriang})
it follows that ${\cal A}(u)$ is invariant by the map $T$, so that
${\cal P}_u$ is a "0\,-1" property.

b) If $\Phi$ is integrable and $\int \Phi \, d\mu \not = 0$, then by
the ergodic theorem $\ad(\Phi)$ reduces to the direction defined by
the mean of $\Phi$. Therefore, when $\Phi$ is integrable, the
interesting case is when $\int \Phi \, d\mu = 0$.

c) If $\Phi$ is a coboundary, $\Phi = \Psi - \Psi \circ T$, then the
set of limit directions can be deduced from the support of the law
of $\Psi$ in $\R^d$. When this law gives a positive measure to each
cone truncated from the origin, then, by ergodicity of $T$, the set
$\ad(\Phi)$ coincide with $\Sp^{d-1}$.

d) Billiards in the plane with periodic obstacles yield geometric
examples of centered vector cocycles with a geometric interpretation
of the limit directions (for these models, see for example
\cite{Pe00}, \cite{SzVa04} for the dispersive billiards,
\cite{Gu10}, \cite{CoGu12} for the billiards with polygonal
obstacles).}
\end{rems}

\begin{lem} There is a closed set $\ad(\Phi)$ such that
$\ad(\Phi)(x) = \ad(\Phi)$, for a.e. $x$. It is empty if and only if
$\Phi$ is a coboundary: $\Phi = \Psi - \Psi \circ T$, with $\Psi$
bounded.
\end{lem} \proof \ Clearly $\ad(\Phi)(x)$ is a closed subset of
$\Sp^{d-1}$. The invariance $\ad(\Phi)(T x) = \ad(\Phi)(x)$ follows
from (\ref{inegTriang}). Using the Hausdorff distance on the set of
closed subsets of $\Sp^{d-1}$ and ergodicity, we obtain that
$\ad(\Phi)(x)$ is a.e. equal to a fixed closed subset.

If $\ad(\Phi)$ is empty, then, for a.e. $x$, the sequence
$(\Phi_n(x))$ is bounded. This implies that there is a measurable
function $\Psi$ such that $\Phi = \Psi - \Psi \circ T$. By
ergodicity of $T$, $\Psi$ is bounded. The converse is clear. \eop

\begin{defi} {\rm $\ad(\Phi)$ will be called {\it set of limit directions} of
$(\Phi_n)$. We write also ${\cal D}(T, \Phi)$ instead of $\ad(\Phi)$
to explicit the dependence on $T$.}
\end{defi}

In other words, the "limit set" $\ad(\Phi)$ is in the transient case
the attractor in the sphere $\Sp^{d-1}$ of the process $(\tilde
\Phi_n)_{n \geq 1}$ introduced above.

In the last section, we will show that under a strong stochastic
hypothesis, this process $(\tilde \Phi_n)_{n \geq 1}$ visits any non
empty open set in $\Sp^{d-1}$ and stays there during longer and
longer intervals of time. This property can be formalized as
follows.

Let $(Z_n)$ be a process defined on $(X, \mu)$ with values in a
metric space $Y$. A first question is about transitivity: does
$(Z_n)$ visit every non empty open set in $Y$. With the previous
notion of limit direction, for the process $({\varphi_n(.) \over
\|\varphi_n(.)\|})$ associated to a transient cocycle $(\varphi_n)$,
this means $\ad(\varphi) = \Sp^{d-1}$.

A stronger quantitative property is the following:
\begin{eqnarray}
&& \limsup_n \frac1n \sum_{1}^n {\bf 1}_V(Z_k(x)) = 1 \text{ a.e.},
\text{ for every non empty open set } V \text { in } Y.
\label{visitV}
\end{eqnarray}
Clearly this property implies $\liminf_n \frac1n \sum_{1}^n {\bf
1}_V(Z_k(x)) = 0$ a.e., for every non empty open subset $V$ in $Y$
with a complement with non empty interior.

We will now discuss some general properties of the set of limit
directions. The set of limit directions $\ad(\Phi_B, T^B)$ for the
induced cocycle $\Phi^B$ and the induced map $T^B$ is denoted by
$\ad(\Phi_B)$ or $\ad(B)$.

We have the equivalence:
\begin{lem} \label{aDCob} a) A cocycle is a coboundary if and only
if there is $B$ of positive measure such that the set of limit
directions  for the induced cocycle on $B$ is empty.

b) If $\Phi$ and $\Phi'$ are cohomologous, there is $B$ such that
the corresponding induced cocycles on $B$ have the same set of limit
directions.
\end{lem} \proof \  a) If $u$ is a limit direction for $(T_A, \Phi_A)$,
then it is also a limit direction for $(T, \Phi)$; hence the
inclusion ${\cal D}(T_A, \Phi^A) \subset {\cal D}(T, \Phi)$.

By a compactness argument on the set of directions, if $(A_n)_{n
\geq 1}$ is a sequence of decreasing sets with positive measure in
$X$, then $({\cal D}(T_{A_n}, \Phi^{A_n}))_{n \geq 1}$ is decreasing
and the intersection is non empty, except if ${\cal D} (T_{A_{n_0}},
\Phi^{A_{n_0}})$ is empty for some $n_0$.

If $(\Phi_n)$ is bounded, or equivalently if $\Phi = \Psi - \Psi
\circ T$, with a bounded $\Psi$, then clearly ${\cal D}(T, \Phi)$ is
empty.

Suppose now that $\Phi$ is a coboundary, $\Phi = \Psi - \Psi \circ
T$, with $\Psi$ measurable. Let $B$ be a set such that $\Psi$ is
bounded on $B$. Then the induced cocycle $(\Phi_n^B)_{n \geq 1}$  is
bounded, since $\Phi_n^B = \Psi - \Psi \circ T_B^n$. Therefore
${\cal D}(T_B, \Phi^B)$ is empty.

Conversely, if there is $B$ of positive measure in $X$ such that
${\cal D}(T_B, \Phi^B)$ is empty, then the induced cocycle
$(\Phi_n^B)$ is bounded, so $\Phi^B$ is a $T_B$ cocycle. By Lemma
\ref{cobInduced} below, $\Phi$ is a coboundary.

b) Let $\Phi$ and $\Phi'$ be such that $\Phi' = T\Psi - \Psi + \Phi$
for a measurable $\Psi$. Let $B$ such that $\Psi$ is bounded on $B$.
Then ${\Phi'}_n^B = T_B^n \Psi - \Psi + \Phi_n^B$ with $T_B^n \Psi -
\Psi$ bounded, which implies that $({\Phi'}_n^B) $ and $(\Phi_n^B)$
have the same limit directions. \eop

\begin{lem} \label{cobInduced} Let $B$ be such that $X = \cup_{k \geq 0}
T^k B$. If $\Phi^B$ is a $T_B$-coboundary, then $\Phi$ is a
$T$-coboundary.
\end{lem} \proof \ For $\mu$-a.a. $y \in X$ there are a unique
$x \in B$ and an integer $k$, $0\le k < R_B(x)$, such that $y = T^k
x$. Suppose that there is $\Psi$ on $B$ such that: $\Phi^B= \Psi -
\Psi \circ T^B$. We define $\zeta$ on $X$ by taking, for $0\le k <
R_B(x)$, $\zeta(T^k x) = \Psi(x) - \Phi_k(x)$.

We have $\Phi(y) = \zeta(y) - \zeta(T y)$. Indeed, for $y = T^k x$,
$0 \le k <R_B(x)-1$, the relation is satisfied by construction; for
$y = T^k x$ with $k=R_B(x) -1$, the relation follows from the
coboundary relation for the induced cocycle. \eop

\vskip 3mm Now let us show that two sets $A_1$ and $A_2$ have always
non disjoint sets of limit directions, unless $\Phi$ is a
coboundary.

\begin{lem} \label{asympSubset} For any two sets $A_1$ and $A_2$,
there is $B_1 \subset A_1$ such that $\ad(B_1) \subset \ad(A_2)$.
\end{lem} \proof \  Let $B_1 \subset A_1$ be such that $\Phi^{A_2}$
(as defined by (\ref{passage})) satisfies $\Phi^{A_2}(x) \leq C$ on
$B_1$, for some constant $C$. If $\ad(B_1)$ is empty, then $\Phi$ is
a coboundary by Lemma \ref{aDCob}. Let $u$ be a limit direction for
$\Phi^{B_1}$.

The cocycle $(\Phi_n^{A_2}(x))$, for a piece of orbit starting and
ending in $A_2$ and for a special sequence of times, can be written
as $(1) + (2) + (3)$ where
\begin{eqnarray*}
(1) &=& \sum_{t=0}^{R_{B_1}(x) -1} \Phi(T^t x), \\
(2) &=& \Phi_{n_k(x_1)}^{B_1}(x_1), {\rm \ with \ } x_1 = T^{R_{B_1}(x)}x,\\
(3) &=& \Phi^{A_2}(x_2), {\rm \ with \ } x_2 =
T_{B_1}^{n_k(x_1)}x_1.
\end{eqnarray*}

The first term (1) corresponds to the path from $A_2$ to $B_1$. The
second term (2) corresponds to visits of the cocycle induced on
$B_1$ to a neighborhood of $u$ (after normalization) with an
arbitrary large norm (such visits exist because $u$ is a limit
direction for the induced cocycle on $B_1$), and the third (3) to
the path from $B_1$ to $A_2$ with a bounded value of the cocycle by
construction.

If we iterate for a long time the induced cocycle (2), the first
term (which is fixed) and the third (which is bounded) are small
compared with the norm of (2). Then $(1) + (2) + (3)$ gives a value
of the induced cocycle on $A_2$ which satisfy the condition that $u$
is a limit direction for $\Phi^{A_2}$. \eop

\vskip 3mm \subsection{\bf Essential limit directions}

\

The observation that the set of limit directions $\ad(\Phi)$ is not
a "cohomological invariant" motivates the following definition.

\begin{defi} \label{asymptDir} {\rm A direction $u \in \Sp^{d-1}$
is called an {\it essential limit direction} for $\Phi$, if, for
every subset $B$ of positive measure, $u$ is a limit direction for
$\Phi^B$. The set of essential limit directions is denoted by
$\ead(\Phi)$.}
\end{defi}

The set $\ead(\Phi)$ can be seen as a "boundary" for $(\Phi_n)$. It
is invariant by cohomology: if $\Phi_1$ and $\Phi_2$ are
cohomologous, then $\ead(\Phi_1) = \ead(\Phi_2)$.

\begin{thm} \label{connecAsymp} 1) $\ead(\Phi)$ is a closed subset of $\Sp^{d-1}$
which is empty if and only if $\Phi$ is a coboundary. For every $B$
of positive measure, $\ead(\Phi^B) = \ead(\Phi)$.

2) If $(\Phi_n)$ is transient and $\Phi$ is bounded, then
$\ad(\Phi)$ is a closed connected non empty subset of $\Sp^{d-1}$.
\end{thm}
\proof \ 1) We have $\ead(\Phi) = \bigcap \ad(\Phi^B)$, where the
intersection is over the family of all measurable subsets of
positive measure.

Clearly, $\ead(\Phi) \subset \ead(\Phi^B)$. Let $A$ with $\mu(A) >
0$. By Lemma \ref{asympSubset}, there is $B_1 \subset B$ such that
every limit direction for $\Phi^{B_1}$ is a limit direction for
$\Phi^A$. If $u$ is in $\ead(\Phi^B)$, then $u$ is in
$\ad(\Phi^{B_1})$, hence in $\ad(\Phi^{A})$. Therefore $\ead(\Phi^B)
\subset \ad(\Phi^{A})$, for all $A$, which implies $\ead(\Phi^B)
\subset \ead(\Phi)$.

2) Let $u_1$ and $u_2$ be two accumulation points of
$\Phi_n(x)/\|\Phi_n(x)\|$ and $\varepsilon>0$. For a.e. $x$, by
transience, for $N$ big enough, we have $\|\Phi_n(x)\| \geq
\varepsilon^{-1}\|\Phi\|_\infty, \forall n \geq N$. By definition,
there exist $n>m>N$ such that
$d(\Phi_m(x)/\|\Phi_m(x)\|,u_1)<\varepsilon$ and
$d(\Phi_n(x)/\|\Phi_n(x)\|,u_2)<\varepsilon$. Moreover, for every
$k$ between $m$ and $n-1$, one has
\begin{eqnarray*}
&&\|{{\Phi_k(x)}\over{\|\Phi_k(x)\|}}-{{\Phi_{k+1}(x)}
\over{\|\Phi_{k+1}(x)\|}}\|\\ &&\leq {{\|\Phi \circ
T^k(x)\|}\over{\|\Phi_{k+1}(x)\|}}+{\|\Phi_k(x)\|}\left|
{{1}\over{\|\Phi_{k+1}(x)\|}}-{{1}\over{\|\Phi_{k}(x)\|}}\right|
\leq 2\varepsilon.
\end{eqnarray*}
Thus, for every $\varepsilon>0$, one has a finite set
$F_\varepsilon$ of points $\Phi_k(x)/\|\Phi_k(x)\|$ on the unit
sphere that can be used to go from  $u_1$ to $u_2$ with jumps of
length smaller than $2\varepsilon$. Now let $\varepsilon$ tend to
zero and consider $F_\infty$ an accumulation point of
$(F_\varepsilon)_{\varepsilon > 0}$ in the set of compact sets of
the sphere equipped with the Hausdorff metric. The set $F_\infty$ is
a connected compact set containing $u_1$ and $u_2$. \eop

\vskip 3mm Let $\widetilde \ev(\Phi)$ be the smallest vector space
of $\R^d$ containing $\ev(\Phi)$. Using Definition \ref{regul0}, we
have:
\begin{thm} For every non coboundary $\Phi$, $\ead(\Phi)$
contains $\Sp(\widetilde \ev(\Phi))$, the sphere at infinity of
$\widetilde \ev(\Phi)$, and is equal to $\Sp(\widetilde \ev(\Phi))$
if $\Phi$ is a regular cocycle.
\end{thm}

{\bf Remarks and questions.} a) A general question is to find the
set of limit directions and the set of essential limit directions of
a given cocycle. What are the possible shapes of these sets ?

b) The rate of growth of the cocycle plays no role in the
"directional process" associated to a cocycle as defined above. This
rate could be taken into account by introducing a scaling in the
notion of limit directions.

c) Let us call "irreducible" a $\R^d$-cocycle which is not
cohomologous to a cocycle with values in a vector subspace of
dimension $< d$. For an irreducible cocycle $\Phi$ what kind of set
$\ad(\Phi)$ can be ? In particular does there exist a recurrent
cocycle $(\Phi_n)$ such that $\ad(\Phi)$ reduces to two antipodal
points.

This question is related to the following remark. Let $\Phi$ and
$\Psi$ with values in $\R^d$ be given. We say that {\it the cocycle
$(\Phi_n)$ dominates $(\Psi_n)$, if there are $C$ and $K$ such that
$\|\Psi_n(x)\| \leq C \|\Phi_n(x)\| +K, \forall n$}.

Clearly this is the case when $\Psi$ is cohomologous to a multiple
of $\Phi$ with a bounded transfer function. The proposition below is
a partial converse.

\begin{prop} Assume that $T_\Phi$ is ergodic on $X \times \R^d$.
If $(\Phi_n)$ dominates $(\Psi_n)$, then $\Psi$ is cohomologous to
$c \, \Phi$ for a constant $c$.
\end{prop} \proof \ Let $I$ be a compact neighborhood of $\{ 0 \}$.
Let $y$ be in $I$. For the times $n_k(x,y)$ such that $y +
\Phi_{n_k(x,y)}(x) \in I$, $|\Psi_n(x)|$ is bounded. The cocycle
$\Psi_n^{Z_I}(x)$ induced of $\Psi$ on the set $Z_I:=X \times I$ is
bounded. Therefore the function $F$ defined on $X \times \R$ by
$F(x,y)= \Psi(x)$ is a coboundary for the map $T_\Phi$: there is
$H(x,y)$ such that $F(x,y)= \Psi(x) = H(T x,y+\Phi(x)) - H(x,y)$.

For every $a \in \R^d$, the function $(x,y)\to H(x, y+a) - H(x,y)$
is $T_\Phi$-invariant, hence a.e. constant by ergodicity of
$T_\Phi$: for a.e. $(x,y)$, there is $c(a)$ such that $H(x, y+a) =
c(a) + H(x,y)$. By the theorem of Fubini, for a.e. $y$,  $H(x, y+a)
= c(a) + H(x,y)$, for almost every $(x, a)$, and $a \to c(a)$ is
Lebesgue measurable. Let us take $y_0$ satisfying this property. We
have $H(x, a+y_0) = c(a) + H(x,y_0), \ \text{for a.e. }(x,a)$;
hence, with $u(a) = c(a-y_0)$ and $h(x) = H(x, y_0)$:
$$H(x, a) = u(a) + h(x), \ \text{for a.e. }(x,a).$$

The relation $H(x, y+a) = c(a) + H(x,y)$ reads: $u(y+a) + h(x) =
c(a) + u(y) + h(x)$ which shows that $u$ is an additive function.

Therefore $H(x,y) = cy + h(x)$ for a constant $c$ and a measurable
function $h$ on $X$ and we have $\Psi(x) = c \,\Phi(x) +h(T x) -
h(x)$. \eop

\vskip 3mm \subsection{\bf A $G_\delta$-property} \label{billmod}

\

Suppose that the map $T = T(\theta)$ and the function defining the
cocycle $\Phi= \Phi^\theta$ depend on a parameter $\theta$. Suppose
that $\Theta$, the set of parameters, is a metric space and that the
dependence of $T(\theta)$ and $\Phi^\theta$ is piecewise continuous.
We denote by ${\cal V}= {\cal V}(u)$ a countable basis of open
neighborhoods of a direction $u$ in $\Sp^{d-1}$.

\begin{thm} Suppose that in the set of parameters
$\Theta$ there is a dense set ${\cal T}$ of values such that the
corresponding set of limit directions is $\Sp^{d-1}$. Then there is
a dense $G_\delta$-set in $\Theta$ with the same property.
\end{thm} \proof \ We can assume that ${\cal T}$ is countable:
${\cal T} = \{\theta_i, i = 1, 2, ... \}$. Fix a direction $u \in
\Sp^{d-1}$. For $\theta \in {\cal T}$, for a.e. $x \in X$, $u$ is a
limit direction for $(\Phi_n^\theta(x))$. Let $K$ be a compact set
of positive measure in $X$ such that for every $i$,
$$ K \subset \{x \in X: u {\rm \ is \
a \ limit \ direction \ for \ } \Phi_n^{\theta_i}(x) \}.$$

For a fixed $x$, for $M \geq 1$ and $V \in {\cal V}(u)$, the set
$$\tilde B_n^{x,V,M} = \{\theta: \|\Phi_n^\theta(x)\| > M  {\rm \
and \ } \Phi_n^\theta(x)/\|\Phi_n^\theta(x)\| \in V\}$$ is an open
set.

If $W$ is an open set in $X$, let
$$\tilde B_{n}^{W, V, M} := \{\theta: \|\Phi_n^\theta(y)\| > M  {\rm \ and \ }
\Phi_n^\theta(y)/\|\Phi_n^\theta(y)\| \in V, \forall y \in W\}.$$

Let $V \in {\cal V}(u)$, $M \in \N$ and $\theta_i \in {\cal T}$. For
each $x \in K$, there exists $n$ such that $\theta_i \in \tilde
B_n^{x, V, M}$. By continuity $\theta_i \in \tilde B_n^{y, V, M}$
for $y$ in an open neighborhood of $x$. Thus there are finitely many
open sets $W_{(V,M)}^1, ..., W_{(V,M)}^{r_i(V,M)}$ covering $K$ and
integers $n_{(V,M)}^1, ..., n_{(V,M)}^{r_i(V,M)}$ such that
$$\theta_i \in \bigcap_{j = 1, ..., r_i(V,M)} \tilde B_{n_{(V,M)}^j}^{W_{(V,M)}^j, V, M}.$$

This proves that, for every $y \in K$, there are $j \in \{1, ...,
r_i(V,M)\}$ and $n_j$ such that
$$\|\Phi_{n_j}^{\theta_i}(y)\| > M  {\rm \ and \ }
\Phi_{n_j}^{\theta_i}(y)/\|\Phi_{n_j}^{\theta_i}(y)\| \in V.$$

For every $i$, $\theta_i$ belongs to the open set $\bigcup_{i}
\bigcap_{j \in \{1, ..., r_i\}} \tilde B_{n_{(V,M)}^j}^{W_{(V,M)}^j,
V, M}$. The dense set of parameters ${\cal T} = \{ \theta_i \}$ is
contained in the countable intersection of open sets:
\begin{eqnarray}
\bigcap_{V \in {\cal V}(u), M \geq 1} \ \bigcup_{i} \bigcap_{j \in
\{1, ..., r_i(V,M)\}} \tilde B_{n_{(V,M)}^j}^{W_{(V,M)}^j, V, M}.
\label{Gdelta}
\end{eqnarray}

Now, suppose that the parameter belongs to the dense $G_\delta$-set
defined above by (\ref{Gdelta}).

For every $V \in {\cal V}(u)$ and $M \in \N$, there is $i$ such that
$$\theta \in \bigcap_{j \in \{1, ..., r_i(V,M)\}} \tilde B_{n_{(V,M)}^j}^{W_{(V,M)}^j,
V, M},$$ i.e., for every $V \in {\cal V}(u)$ and $M \in \N$, there
are $i$ and $j \in \{1, ..., r_i(V,M)\}$ such that
$$\forall y \in W_{(V,M)}^j, \ \|\Phi_{n_{(V,M)}^j}^{\theta}(y)\| > M {\rm \ and \ }
\Phi_{n_{(V,M)}^j}^\theta(y)/\|\Phi_{n_{(V,M)}^j}^\theta(y)\| \in
V.$$

As $W_{(V,M)}^j, j = 1, ..., r_i(V,M),$ is a covering of $K$, for
each point $y \in K$, for all $V \in {\cal V}(u)$, all $M \in \N$,
there is $n$ such that $\|\Phi_{n}^{\theta}(y)\| > M {\rm \ and \ }
\Phi_{n}^\theta(y)/ \|\Phi_{n}^\theta(y)\| \in V$.

Therefore for each $y \in K$, $u$ is a limit direction. As it is a
0\,-1-property, the property that $u$ is a limit direction holds for
a.e. $x$. \eop

\vskip 3mm \subsection{\bf Limit directions and limit distributions}
\label{limdistr}

\

\begin{lem} \label{limDistrib} Suppose that, for a sequence of integers $(k_n)$
and a sequence $(a_n)$ tending to $\infty$, $(\Phi_{n})$ satisfies a
limit theorem in distribution:
$$a_n^{-1} \Phi_{k_n} \overset{{distrib}} \longrightarrow {\cal L},$$
where ${\cal L}$ is a probability measure on $\R^d$ giving a
positive probability to each non empty open set. Then the set $\Cal
{D}(\Phi)$ of limit directions of $(\Phi_n)$ is $\Sp^{d-1}$. This
applies in particular if  $(\Phi_{n})$ satisfies a non degenerated
CLT for a subsequence and an adapted normalization.
\end{lem} \proof \ Suppose that $u \in \Sp^{d-1}$ is not a limit
direction for $(\Phi_{n})$. By the dichotomy (cf. Remark
\ref{remderiv} a)), there is $M$ and an open regular neighborhood $V
= V(u)$ of $u$ in $\R^d$ such that, a.e. $x$ belongs for some $N$ to
the set
$$C_N:=\{x: \forall n \geq N, \|\Phi_{k_n}(x)\| \leq M {\rm \ or \ }
 \Phi_{k_n}(x) / \|\Phi_{k_n}(x)\| \not \in V(u)\}.$$

The sequence of sets $(C_N)$ is increasing and $\mu(\bigcup_N C_N) =
1$.

From the assumption, we have
$$\liminf_{n \to + \infty}\mu\{x \in X: a_n^{-1} \Phi_{k_n}(x) \in V \}
\geq {\cal L}(V) > 0.$$
Therefore for any $\alpha >0$, there is $N$ such that, for $n \geq
N$, there is a set $B$ in $X$ of measure $ > {{\cal L}(V) \over 2}$
such that: $\|\Phi_{k_n}(x) - a_n u\| \leq \alpha a_n$, which
implies
$$(1- \alpha) a_n \leq \|\Phi_{k_n}(x)\| \leq (1+\alpha) a_n, {\rm \ for \ } x \in B.$$
Hence: $\|\Phi_{k_n}(x)/\|\Phi_{k_n}(x)\| - u \| \leq {\alpha \over
1- \alpha}$ on a set of measure ${{\cal L}(V) \over 2} > 0$. If we
take $N$ such that $\mu(C_N) > 1 - {{\cal L}(V) \over 2}$, there is
a contradiction for $n > N$ big enough.

This applies when ${\cal L} = {\cal N}(0, \Gamma)$ where $\Gamma$ is
a non degenerated covariance matrix.\eop

\vskip 3mm \subsection{\bf Oscillations of 1-dimensional cocycles}
\label{oscildim1}

\

We discuss now the notion of limit directions in the special case of
cocycles with values in $\RR$. About oscillations of 1-dimensional
cocycles, let us mention the work of Derriennic (\cite{De10}) where
other references on the subject, in particular of Tanny, Kesten,
Wos, can also be found. For completeness we give below a short
presentation, related to the notion of limit direction, of some
results on 1-dimensional cocycles.

We consider a conservative transformation $T$ of a space $(X, \mu)$
where $\mu$ is $\sigma$-finite and non singular for $T$. Notice that
in Lemmas \ref{decompLem} and \ref{kest2} below, ergodicity is not
assumed. Recall that equalities between functions are understood
$\mu$-a.e.

In this subsection, we denote by $\varphi$ and $(\varphi_n)$
respectively a given measurable function and the corresponding
cocycle over $(X, \mu, T)$. We define two sets, clearly
$T$-invariant (cf. (\ref{inegTriang})):
\begin{eqnarray}
F_\varphi^+ := (\inf_n \varphi_n(.) > -\infty), \ F_\varphi^- :=
(\sup_n \varphi_n(.) < +\infty). \label{defFpm}
\end{eqnarray}

Let us recall the following classical lemma (filling scheme).
\begin{lem} \label{decompLem} There are two functions $h^+$ and
$g^+$ defined on $F_\varphi^+$ (resp. $h^-$ and $g^-$ defined on
$F_\varphi^-$) with values in $[0,+\infty[$ such that
\begin{eqnarray}
\varphi(x) &=& h^+(T x)- h^+(x) + g^+(x), \for \mu \ae x \in
F_\varphi^+, \label{decomp}\\
\varphi(x) &=& -h^-(T x)+ h^-(x) - g^+(x), \for \mu \ae x \in
F_\varphi^-, \label{decomp2}
\end{eqnarray}
On the invariant set $F_\varphi^{+\infty} := (\varphi_n(.)
\rightarrow +\infty)$, we have $\sum_{k= 0}^\infty g^+(T^k x)=
+\infty$.

If $(\varphi_n)$ is recurrent, then the space $X$ decomposes in two
invariant sets, each of them possibly of zero measure, one on which
$\varphi$ is a coboundary, the other on which $\sup_n \varphi_n(.) =
+\infty$ and $\inf_n \varphi_n(.) = -\infty$.
\end{lem} \proof \ Let $m_n(x) := \min_{1 \leq k \leq n}
(\varphi_k(x))$, $n \geq 1$. We have
\begin{eqnarray*}
m_{n+1}(x)  = \min(\varphi(x), \varphi(x) + m_n(T x)) = \left\{
\begin{array}{ll} \varphi(x) - m_n^-(T x), \ifrm m_n(T x) \leq
0, \\  \varphi(x)= \varphi(x) - m_n^-(T x), \ifrm m_n(T x) > 0,
\end{array} \right.\end{eqnarray*} which implies $m_{n+1}(x) =
\varphi(x) - m_n^-(T x)$. Since the limit $m_\infty(x) := \lim_n
m_n(x)$ is finite on $F_\varphi^+$, it follows:
\begin{eqnarray*}
\varphi(x) =  m_{\infty}^-(T x) - m_{\infty}^-(x) + m_\infty^+(x), \
x \in F_\varphi^+.
\end{eqnarray*}

This gives the decomposition (\ref{decomp}) on $F_\varphi^+$, with
$h^+ = m_{\infty}^-$ and $g^+ = m_\infty^+$. We get (\ref{decomp2})
by changing $\varphi$ into $-\varphi$.

On the invariant set $F_\varphi^{+\infty} = (\varphi_n(.)
\rightarrow +\infty)$ the decomposition (\ref{decomp}) holds. Let us
show that $\sum_{k= 0}^\infty g^+(T^k x)= +\infty$.

We have $\varphi_n(x) = h^+(T^n x) - h^+(x) + \sum_{k= 0}^{n - 1}
g^+(T^k x)$. Let $M_K$ be the subset of $F_\varphi^{+\infty}$ where
$h^+$ is bounded by a finite constant $K$. By conservativity of $T$,
for a.e. $x$ in $M_K$ there is a subsequence $(n_j(x))$ such that
$T^{n_j(x)}(x) \in M_K$ and therefore $\sum_{k= 0}^{n_j(x) - 1}
g^+(T^k x) \to +\infty$. It follows $\sum_{k= 0}^{\infty} g^+(T^k x)
= +\infty$ for a.e $x$ in $M_K$ and, since $K$ is arbitrary,
$\sum_{k= 0}^{\infty} g^+(T^k x) = +\infty$ a.e.

Suppose now that $(\varphi_n)$ is recurrent. Let $g_\infty^+(x) =
\sum_0^\infty g^+(T^k x) \in [0, +\infty]$. The induced cocycle
$(\varphi_n^B)$ is recurrent for any set $B$ on which $h^+$ is
bounded and therefore $g_\infty^+(x) < +\infty$, a.e. on $B$. As the
sets $B$ cover $F_\varphi^+$, this implies $g_\infty^+(x) <
+\infty$, a.e. on $F_\varphi^+$. Therefore $g^+(x) = g_\infty^+(x) -
g_\infty^+(T x)$, and the restriction of $\varphi$ to the invariant
set $F_\varphi^+$ is a coboundary. Likewise, $\varphi$ is a
coboundary on the invariant set $F_\varphi^-$.

So we have proved that, for any recurrent cocycle, the space $X$
decomposes in two sets, the set $F_\varphi^+ \cup F_\varphi^-$ on
which $\varphi$ is a coboundary and its complement on which
$\varphi_n$ oscillates between $+\infty$ and $-\infty$. \eop

\vskip 3mm The lemma implies that $\varphi$ is a coboundary on the
invariant set $\{ x: \varphi_n(x) {\rm \  is \ bounded}\}$. It is
well known that, if $(\varphi_n)$ is uniformly bounded, then the
transfer function is bounded.

The previous lemma gives a simple way to prove and to slightly
extend a result of Kesten on the rate of divergence in dimension 1
of a non recurrent cocycle. We consider a conservative dynamical
system with a $\sigma$-finite invariant measure $\mu$. The
$\sigma$-algebra of $T$-invariant sets is denoted by ${\mathcal I}$.

As $\mu$ is $\sigma$-finite, we can choose a function $p$ on $X$
such that $\mu(p) = 1$ and $0 < p(x) \leq 1, \forall x$. By the
ratio ergodic theorem $\lim_{n \to\infty}
\frac{\sum_{k=0}^{n-1}f(T^kx)} {\sum_{k=0}^{n-1}p(T^kx)}
=\EE_{p\mu}[{f \over p} |\mathcal{I}](x)$, $\mu$-a.e. $x\in X$, for
$f \in \mathbb{L}^1(\mu)$. We have $p_n(x) = \sum_0^{n-1} p(T^k x)
\to \infty$ since the system is conservative.

\begin{lem} \label{kest2} (cf. also \cite{Ke75}, \cite{De10})
Suppose that the  $T$-invariant measure $\mu$ is conservative
$\sigma$-finite.

1) If $\varphi$ is a nonnegative measurable function, then for
$\mu$-a.e. $x$ the sum $\sum_0^\infty \varphi(T^k x)$ is either 0 or
$+\infty$, and $\liminf_n {\varphi_n(x) \over p_n(x)}\in  \, ]0,
+\infty]$ on the set $\{\sum_0^\infty \varphi(T^k .) \not = 0\}$.

2) For any measurable $\varphi$, $\liminf_n {\varphi_n(x) \over
p_n(x)} > 0$ for $\mu$-a.e. $x$ on the set $(\varphi_n(.)
\rightarrow +\infty)$.

The cocycle $(\varphi_n)$ is recurrent on the invariant set $\{x:
\lim{\varphi_n(x) \over p_n(x)} = 0\}$.

3) If $\mu$ is a $T$-invariant probability measure, then for any
measurable $\varphi$, $\liminf_n {1\over n} \varphi_n(x)
> 0$ for $\mu$-a.e. $x$, on the set $(\varphi_n(.)
\rightarrow +\infty)$ and the cocycle $(\varphi_n)$ is recurrent on
the invariant set $\{x: \lim{1\over n} \varphi_n(x) = 0\}$.
\end{lem} \proof \ 1a) First, suppose that $\varphi=1_B$ where $B$
is a measurable subset. For $N \geq 0$, let $B^{(N)} := \cup_0^N
T^{-k} B$, $B^{(\infty)} = \cup_0^\infty T^{-k} B$. We have $T^{-1}
B^{(\infty)} \subset B^{(\infty)}$, hence $T^{-1} B^{(\infty)} =
B^{(\infty)}$ up to a negligible set as $(X,\mu, T)$ is
conservative.

On the complement of $B^{(\infty)}$, we have $\sum_0^n 1_B(T^k x) =
0$. By the ratio ergodic theorem $\lim_n{1\over p_n(x)} \sum_0^{n-1}
1_{B}(T^k x) \, p(T^k x) = \EE_{p\mu}(1_{B}|{\mathcal {\mathcal I}})
(x)$.

As $\sum_{j=0}^{L-1} 1_B(T^j x) \geq 1_{B^{(L)}}(x)$ and $\lim_n
p_n(x) = +\infty$, we have
\begin{eqnarray*}
&&\liminf_n {1 \over p_n(x)} \sum_{k=0}^{n-1} 1_B(T^k x) \geq \frac
1L \liminf_n {1 \over p_n(x)} \sum_{k=0}^{n-1} 1_{B^{(L)}}(T^k x)
\\&& \geq \frac 1L \lim_n {1 \over p_n(x)} \sum_{k=0}^{n-1}
1_{B^{(L)}}(T^k x) \, p(T^k x) =
 \frac 1L \EE_{p\mu}(1_{B^{(L)}}|{\mathcal I})(x).
\end{eqnarray*}
Therefore from the relation
\begin{eqnarray*}
\bigcup_L \uparrow \{x: \EE_{p\mu}(1_{B^{(L)}}|{\mathcal I})(x) > 0
\} = \{x: \EE_{p\mu}(1_{B^{(\infty)}}|{\mathcal I})(x) > 0 \} =
1_{B^{(\infty)}},
\end{eqnarray*}
it follows
\begin{eqnarray}
\liminf_n {1 \over p_n(x)} \sum_{k=0}^{n-1} 1_B(T^k x) > 0 {\rm \ on
\ } B^{(\infty)}. \label{posi1}
\end{eqnarray}

1b) Now let $\varphi$ be any nonnegative function. For $j \in \Z$,
let $B_j:=\{2^j \leq \varphi < 2^{j+1}\}$. We get: $\sum_0^\infty
\varphi(T^k x) > 0 \Leftrightarrow \exists j: x \in B_j^{\infty}$.
Using (\ref{posi1}) applied to the sets $B_j$ and the inequality
$\varphi \geq \sum_{j= -\infty}^{+\infty} 2^j 1_{B_j}$, we obtain:
$$\liminf_n {1\over p_n(x)} \sum_0^{n-1} \varphi(T^k x) \geq \sum_{j=
-\infty}^{+\infty} 2^j \, \liminf_n {1\over p_n(x)} \sum_0^{n-1}
1_{B_j}(T^k x) > 0, {\rm \ on  \ } (\sum_0^\infty \varphi(T^k .)
> 0).$$

2) As (\ref{decomp}) of Lemma \ref{decompLem} holds on the set
$F_\varphi^{+\infty} = (\varphi_n(.) \rightarrow +\infty)$ and
$\sum_{k= 0}^\infty g^+(T^k x) = +\infty$, we can apply 1) to $g^+$.
Since $\varphi_n(x) \geq - h^+(x) + g_n^+(x)$, we get:
$$\liminf { \varphi_n(x) \over p_n(x)} \geq \liminf_n {g^+_n( x) \over p_n(x)}
> 0, \ \for \mu \ae \in F_\varphi^{+\infty}.$$

This implies that $(\varphi_n(.))$ is recurrent on the invariant set
$\{x: \lim{\varphi_n(x)\over p_n(x)} = 0 \}$, since we can not have
$\varphi_n(x) \rightarrow +\infty$ or $\varphi_n(x) \rightarrow
-\infty$ on this set.

In particular, if $T$ is ergodic and $\varphi$ integrable with
$\mu(\varphi) = 0$, we have $\lim_n {\varphi_n(x)\over p_n(x)} =
\mu(\varphi) = 0$ which implies the recurrence of the cocycle
$(\varphi_n)$.

3) When the measure is finite, then we take $p(x) =1$ and $p_n(x)$
is replaced by $n$.\eop

\vskip 3mm \begin{proposition} For a 1-dimensional cocycle
$(\varphi_n)$ generated by $\varphi$ over an {\it ergodic} dynamical
system, if $\varphi$ is not a coboundary, one of the following
(exclusive) properties is satisfied: \hfill \break 1)
$\limsup_n\varphi_n = -\liminf \varphi_n = + \infty$, \hfill \break
2) $\varphi_n$ tends to $+\infty$, \hfill \break 3) $\varphi_n$
tends to $-\infty$.
\end{proposition}
\proof \ By ergodicity, with the notation (\ref{defFpm}), one of the
sets  $(F_\varphi^+ \cup F_\varphi^-)^c$, $F_\varphi^+$,
$F_\varphi^-$ has full measure.

The first case is equivalent to property 1) and to the equality
$\ead(\varphi) = \{-\infty, +\infty\}$. Suppose now that
$F_\varphi^+$ has full measure. Then we have the decomposition
(\ref{decomp}), Lemma \ref{decompLem}, with equality a.e. Hence
$\varphi_n(x) \geq - h^+(x) + g_n^+(x), \for \mu \ae x $. Since
$\varphi$ is not a coboundary, $g^+$ is non negligible. This implies
property 2). Likewise 3) holds if $F_\varphi^-$ has full measure.
\eop

\vskip 3mm This leads to the following remarks. \hfill \break - if
the cocycle $(\varphi_n)$ is recurrent, it oscillates between
$+\infty$ and $-\infty$, unless $\varphi$ is a coboundary with a
transfer function bounded from above or below; \hfill \break - if
$\varphi$ is a coboundary, $\varphi = T \psi - \psi$, then if $\psi$
is not essentially bounded from above (resp. from below), then
$\limsup_n\varphi_n = +\infty$ (resp. $\liminf \varphi_n = -
\infty$); \hfill \break - we have $- \infty \not \in \Cal D
(\varphi)$ if and only if $\varphi = Th - h + g$, with $h$, $g$ non
negative and $g$ non negligible. This is equivalent to $\lim_n
\varphi_n = + \infty$; \hfill \break - $ \Cal D (\varphi) $ is empty
if $\varphi$ is a coboundary, $\varphi = Th - h $, with $h$
essentially bounded.

 For the constructions below, we need the following lemma:
\begin{lem} \label{over} Let $(\ell_n)$ be a strictly increasing sequence of
integers. For any ergodic dynamical system there is $h$ non negative
such that, for a.e. $x$, $h(T^n x) \geq \ell_n$ infinitely often.
\end{lem}
\proof \ There exists a strictly increasing sequence of positive
real numbers $(c_j)$ and a strictly increasing sequence of natural
numbers $(n_j)$, both tending to infinity, and a non negative
measurable function $f$ such that
$$\lim_j {1\over c_j n_j} \sum_{k=1}^{n_j} f(T^k x) = +\infty, \text{
a.e. }$$ We put $d_k = c_j$, for $n_{j-1} \leq k < n_{j}$. For a.e.
$x$, for $j$ big enough, we can define a non decreasing sequence
$(k_j(x))$ such that $\lim_j k_j(x) = + \infty$ and $f(T^{k_j(x)}x)
\geq c_j \geq d_{k_j(x)}$. (Put $k_j(x) := \max \{k \leq n_j:
f(T^{k}x) \geq c_j\}$.)

Now we can define a non decreasing function $\gamma$ on $\R^+$ by
putting $\gamma(y) = \ell_k$, for $d_k \leq y < d_{k+1}$. In
particular, we have $\gamma(d_k) = \ell_k$.

Let $h(x) := \gamma(f(x))$. Then we have: $h(T^k x) = \gamma(f(T^k
x)) \geq \gamma(d_k) = \ell_k$, if $f(T^k x) \geq d_k$.

Therefore for a.e. $x$, for $j$ big enough, $h(T^{k_j(x)}(x)) \geq
\ell_{k_j(x)}$. \eop

\vskip 3mm For every $B \subset X$ of positive measure, $\varphi$ is
a coboundary, if and only if the induced cocycle $\varphi^B$ is a
coboundary for the induced map $T_B$. If $\varphi$ is not a
coboundary, then the inclusion $\ad(\varphi^B) \subset \ad(\varphi)$
is general, but can be strict.

{\it Example:} Let $B$ with $\mu(B) > 0$. Take $\varphi = Th - h +
1$, with $h(x) \geq 0$ on $B$. Then, if $R_n(x)$ denotes the $n$-th
return time to $B$, we have
$$\varphi_n^B(x) = h(T_B^n x) - h(x) + R_n(x) \geq - h(x) + R_n(x) \to +\infty.$$
We can choose $h$ on $B^c$ such that, for the cocycle $\varphi_n(x)
= h(T^n x) - h(x) + n$, we have $h(T^n x) \leq -n^2$ infinitely
often (Lemma \ref{over}). Therefore $-\infty \in \ad(\varphi) \not =
\ad(\varphi_B) = \{+\infty\}$.

\vskip 3mm {\it Reverse cocycle}

Recall that the reverse cocycle $(\check \varphi_n)_{n\geq 0}$ is
defined by $\check \varphi_0 = 0$ and
$$\check \varphi_n(x) = -\varphi_n(T^{-n} x) = - \sum_{k= 1}^{n}
\varphi(T^{-k} x), \  \text{ for } n \geq 1.$$

 For an ergodic system, if $\varphi$ is integrable and if $\lim_n
\varphi_n = + \infty$, then $\lim_n \check \varphi_n  = - \infty$,
since both conditions are equivalent to $\mu(\varphi) > 0$.

If $\varphi$ is non integrable, we can have  $\lim_n {\varphi_n} = +
\infty$ and $\limsup_n \check \varphi_n  = + \infty$. (see also
\cite{De10}).

{\it Example:} Let $\varphi =Th - h + 1$, with $h$ non negative. We
have $\lim_n \varphi_n(x) = +\infty$. The reverse cocycle reads
$$\check \varphi_n(x) = -\varphi_n(T^{-n}x) = -h(x) + h(T^{-n}x) - n.$$
If $h$ is chosen such that the inequality $h(T^{-n} x) \geq n^2$
occurs infinitely often for a.e. $x$ (Lemma \ref{over}), then
$+\infty$ is a limit direction for the reverse cocycle.

\vskip 3mm {\bf A result in dimension 2}

Let us mention a partial result for 2-dimensional cocycles:
\begin{prop} \label{asympdirect} Let $\Phi: X \to \RR^2$ be an
integrable and centered function. If $(\Phi_n)$ is a transient
cocycle over an ergodic dynamical system, then $\ad(\Phi) \cup
(-\ad(\Phi)) = \Sp^1$ for a.e. $x\in X$.
\end{prop} \proof \ We denote by $<u,v>$ the scalar product in $\RR^2$.
Let $v\in \Sp^1$, and let $v^{\perp}\in \Sp^1$ be such that $\langle
v,v^{\perp}\rangle = 0$. The function $x \mapsto
\langle\Phi(x),v^{\perp}\rangle$ has zero integral and $T$ is
ergodic, hence the cocycle $\langle(\Phi_n),v^{\perp}\rangle$ is
recurrent. For a.e. $x\in X$ there is a sequence $n_k(x) \to \infty$
and $c>0$ such that $|\langle \Phi_{n_k}(x),v^{\perp} \rangle|<c$.
As $(\Phi_n)$ is transient, for a.e. $x$ $|\langle \Phi_{n_k}(x),v
\rangle|$ is not bounded. There is a subsequence $(n_{k_j})$ such
that $\Phi_{n_{k_j}}(x)/\|\Phi_{n_{k_j}}(x)\|$ converges to $v$ or
$-v$, i.e., $v$ or $-v \in \ad(\Phi)(x)$. \eop

By what precedes and Theorem \ref{connecAsymp}, when $\Phi$ is
bounded, the set $\ad(\Phi)$ is an arc of a circle with length $\geq
\frac 12$.

\section{\bf Application of the CLT, martingales,
invariance principle \label{applCLT}}

For a large class of dynamical systems of hyperbolic type, the
method introduced by M. Gordin in 1969 gives a way to reduce, up to
a regular coboundary, a Hölderian function $\Phi$ to a function
satisfying a martingale condition. This allows to prove for regular
functions which are not coboundaries, not only a CLT, but also a CLT
for subsequences of positive density and the functional CLT (or the
invariance principle).

In this subsection, we recall some results for martingale increments
and briefly mention their application to find the set of essential
directions..

 \vskip 3mm \subsection{\bf Martingale methods and essential limit directions}
\label{chaotMap}

\

The theorem of Ibragimov and Billingsley stated in terms of
dynamical systems, gives a CLT which can be extended to several
improvements:

\begin{proposition} \label{CLTMart} Let $(X,{\cal A},\mu,T)$ be an ergodic
invertible dynamical system and ${\cal F}$ a sub $\sigma$-algebra of
${\cal A}$ such that $ {\cal F} \subset T^{-1} {\cal F}$. Let $\Phi$
be a $\R^d$-valued square integrable function, ${\cal F}$-measurable
and such that the sequence $(\Phi\circ T^n)_{n \in \Z}$ is a
sequence of martingale increments with respect to $(T^{-n} {\cal
F})$ (equivalently by stationarity: $\E(\Phi |T {\cal F}) = 0$).

If $\Phi$ is non contained a.s. in a fixed hyperplane, the cocycle
$(\Phi_n)$ is such that  $({1\over \sqrt {n}} \Phi_n)_{n \geq 1}$
has asymptotically a Gaussian law, with a non degenerated covariance
matrix $\Gamma$.

 For every strictly increasing sequence of measurable functions
$(k_n)_{n \ge 1}$ with values in $\N$ such that, for a constant $a
\in ]0, \infty[$, $\lim_n { k_n(x) \over n} = a$ exists a.e. we
have:
$$ {1\over \sqrt {n}} \Phi_{k_n(.)}(.) \ {\buildrel {\cal L} \over {\rightarrow}}
\ {\cal N}(0, a^{-1} \Gamma).$$ Moreover the cocycle $(\Phi_n)$
satisfies the invariance principle.
\end{proposition}

\vskip 3mm
\begin{thm} \label{essasympCLT} $\ead(\Phi) = \Sp^{d-1}$ under the conditions
of the previous proposition.
\end{thm}
\proof \ Let $B$ be a subset of positive measure and let $(R_n(x))$
be the sequence of visit times in $B$. The induced cocycle
$(\Phi_n^B)$ is obtained by sampling the cocycle $(\Phi_n)$ at the
random times $R_n$ of visits to $B$.

We have (Kac lemma): $\lim_n {R_n(x) \over n}= {1 \over \mu(B)}$,
since by the ergodic theorem:
$${n \over R_n(x)} = {1 \over R_n(x)} \sum_{j=0}^{R_n(x) - 1}
1_B(T^j x) \to \mu(B).$$ Therefore $(\Phi_n^B)$ satisfies the CLT,
with the same covariance matrix as for the cocycle $(\Phi_n)$ up to
a scalar. By Lemma \ref{limDistrib}, this implies the result. \eop

If the covariance matrix is degenerated, the set $\ead(\Phi)$ is the
unit sphere of a subgroup isomorphic to $\R^{d'}$, for $d' < d$.

\vskip 3mm {\it Reduction by cohomology to martingale increments}

When Gordin's method can be applied, using Theorem \ref{essasympCLT}
and the fact that the set of essential limit directions is the same
for two cohomologous cocycles, we obtain $\ead(\Phi) = \Sp^{d'-1}$,
for some $d' \leq d$, if $\Phi$ is  Hölderian with values in $\R^d$.

This method can be used for Hölderian functions in many systems,
among which: piecewise continuous expansive maps of the interval,
toral automorphisms, geodesic and diagonal flows on homogeneous
spaces of finite volume, dispersive billiards in the plane.

Let us give an explicit example.

\begin{proposition} Let $T$ be an ergodic endomorphism of the torus $\T^r$, $r\geq 1$,
endowed with the Lebesgue measure. If $\Phi$ is a Hölderian function
with values in $\R^d$, then $\ead(\phi) = \Sp^{d'-1}$, for $d' \leq
d$.
\end{proposition}
\proof \ The function $\Phi$ is cohomologous to $\Psi$ such that
$(\Psi \circ T^n)$ is a sequence of vector $d'$-dimensional
martingale increments (with $d' \leq d$) (cf. \cite{LB99}). We can
apply Proposition \ref{CLTMart}, then Theorem \ref{essasympCLT}.
\eop

The situation for the models where Gordin's method is available is
comparable to that of cocycles which are regular in the sense of
Definition \ref{regul0}.

In the last section we will deduce a stronger property from the
invariance principle.

\vskip 3mm \goodbreak \subsection{\bf Invariance principle and
behavior of the directional process} \label{invprincSect}

\

Let $(X, \mu, T)$ be an ergodic dynamical system and let $\Phi$ be a
measurable function on $X$ with values in $\mathbb{R}^d$, $d \geq
2$. In this subsection, $\Phi$ is assumed to be bounded and
centered.

We are going to give conditions on $(\Phi_n)$ which imply the
property (\ref{visitV}) introduced in the first section for the
directional process $Z_n = \tilde \Phi_n(.) = {\Phi_n(.) \over
\|\Phi_n(.)\|}$.

We denote by $(W_n^\Phi)_{n \geq 1}$ (or simply $(W_n)_{n \geq 1}$)
the interpolated piecewise affine process with continuous paths
defined for $x \in X$ and $n \geq 1$ by
$$W_{n}^\Phi(x,s)=\Phi_k(x)+ (ns-k) \Phi(T^k x)
\ {\rm if} \ s\in [{k \over n},{k+1 \over n}[.$$

If $C$ is a cone with non empty interior $\overset{\circ} C$ and
boundary $\partial C$ of measure 0, the amount of time spent by
$W_n^\Phi(x,s)$ in $C$ is
\begin{eqnarray*}
\tau_{n,C}^\Phi(x)=\int_0^1{\bf 1}_C(W_n^\Phi(x,s)) \, ds.
\end{eqnarray*}

Recall that the {\it invariance principle} for $\Phi$, sometimes
called Donsker's invariance principle, means here that the process
$({W_n^\Phi(x,.) \over \sqrt{n}})_{n \geq 1}$ (defined on the
probability space $(X,\mu)$ and with values in the space $(\Cal
C_d([0,1], \| \ \|_\infty)$ of continuous functions from $[0,1]$ to
$\R^d$ endowed with the uniform norm) converges in distribution to
the standard Brownian motion in $\R^d$ (cf. \cite{Bi}).

As mentioned before, the invariance principle, a by-product of the
martingale method, is valid for large classes of regular functions
in many dynamical systems of hyperbolic type.

\vskip 3mm
\begin{thm} \label{invPrincThm} Suppose that $(X, \mu, T)$ is
ergodic, that the invariance principle is satisfied for a centered
bounded function $\Phi\ :\ X \rightarrow \R^d$ and that $C$ is a
cone with non empty interior in $\R^d$ and with a complement with
non empty interior. Then, for almost every $x$,
$$\limsup_{n \rightarrow \infty} \tau_{n,C}^\Phi(x)=1\ {\rm and}\
\liminf_{n \rightarrow \infty} \tau_{n,C}^\Phi(x)=0$$
\end{thm}

We need preliminary lemmas before the proof of Theorem
\ref{invPrincThm}. Firstly, let us remark that the property stated
in the theorem holds for the Brownian motion.

\vskip 3mm {\bf Visit of the Brownian motion in cones}

Let $(B_s)$ denote the standard Brownian motion in $\R^d$, for $d
\geq 2$. Consider a cone $C$ with non empty interior
$\overset{\circ} C$ in $\R^d$ and with a complement with non empty
interior. The amount of time spent by $(B_s)$ in $C$ during the
interval $[0, t]$ is $\displaystyle \tau_{C}(t) = \int_0^t{\bf
1}_C(B_s) \, ds$.

\begin{prop} \label{BrownCase} We have a.s. $\limsup_{t\rightarrow \infty}
{\tau_C(t) \over t}=1$ and $\liminf_{t\rightarrow \infty} {\tau_C(t)
\over t}=0$.
\end{prop}
\proof \ \ Since the variable $\limsup_{t\rightarrow \infty}
{\tau_C(t) \over t}$ is asymptotic, it is  a.s. equal to a constant
value $\ell \in [0,1]$. Because of the scaling property of the
Brownian motion and because $C$ is a cone, we have
\begin{eqnarray*}
\PP({\tau_C(t) \over t}\in I)&=&\PP([\int_0^t{\bf 1}_C(B_s) \, {ds
\over t}] \in I) = \PP([\int_0^1{\bf 1}_C(B_{ts}) \, ds] \in I)\\
&=&\PP([\int_0^1{\bf 1}_C(\sqrt{t}B_{s}) \, ds] \in I) =
\PP([\int_0^1{\bf 1}_C(B_{s}) \, ds] \in I).
\end{eqnarray*}

Take $\alpha \in (0,1)$. As the cone $C$ has a non empty interior,
we have $\PP(B_\alpha \in \overset{\circ} C) > 0$ and, knowing that
$B_\alpha$ is in $\overset{\circ} C$, we also have
$\PP(B_s\in\overset{\circ} C,\ \forall s\in(\alpha,1)) > 0$. The
obvious inequality
$$\PP({\tau_C(t) \over t} > 1-\alpha)\geq\PP(B_s\in\overset{\circ} C,\
\forall s\in(\alpha,1))$$
then implies
\begin{eqnarray}
\PP({\tau_C(t) \over t} > 1-\alpha)>0, \, \forall \alpha>0.
\label{Palpha}
\end{eqnarray}

We have $\limsup_t \PP({\tau_C(t) \over t}> \ell+\varepsilon)\leq
\PP( \limsup ({\tau_C(t) \over t}> \ell+\varepsilon))=0, \forall
\varepsilon
> 0$.

Now the distribution of ${\tau_{C}(t)\over t}$ does not depend on
$t$, so that we have $\PP({\tau_{C}(t)\over t}>\ell+\varepsilon)=
0$. In view of (\ref{Palpha}), this implies that $\ell+\varepsilon >
1-\alpha$. But $\alpha$ and $\varepsilon$ being arbitrary small, one
gets  $\ell\geq 1$, that is $\ell=1$. By considering the complement,
we obtain the result for  $\liminf$. \eop

This suggests that, if we can approximate our process $(W_n^\Phi)$
by a Brownian motion, then the property claimed in Theorem
\ref{invPrincThm} holds. This is the case, for example if we can
assert that for every $\gamma > 1/4$, there exists $C>0$,  so that,
for all $t\in[0,1]$, one has a.s.
$$\|B(nt)-W_n^\Phi(t)\|\leq Cn^\gamma.$$
Such a property is sometimes called an almost sure invariance
principle. It has been established for some hyperbolic or
quasi-hyperbolic systems (see Gouëzel \cite{Go10}). To deduce the
desired property for $W_n$ from the one satisfied by the Brownian
motion, we need to control the amount of time spent by the Brownian
motion not too far from the origin and to enlarge or shrink the cone
we are interested in to get convenient estimates. We will not do
these computations here because they are very similar to what is
done below.

Indeed, we will show that the "plain" Donsker's invariance principle
suffices. From the preceding proof for the Brownian motion, we just
keep in mind that (\ref{Palpha}) in Proposition \ref{BrownCase} is
true.

We need to know that, most of the time, $W_n$ is far from the
origin:

\begin{lem} \label{retcompact} If $\Phi$ is not a coboundary,
for every $M > 0$, the asymptotic frequency of visits of the process
$(W_n)_{n \geq 1}$ to the ball $B(0, M)$ with center at the origin
and radius $M > 0$ in $\R^d$ is almost surely zero:
\begin{eqnarray}
\lim_n \int_0^1 {\bf 1}_{B(0,M)}(W_n(x,s)) \, ds = 0, \text{ for
a.e. } x. \label{asympfreq}
\end{eqnarray}
\end{lem}
\proof \ For $K >0$, the ergodic theorem applied to $(X \times \R^d,
T_\Phi, \lambda = \mu \times dy)$ and ${\bf 1}_{B(0, K)}$ ensures
the existence for a.e. $(x,y) \in X \times \R^d$ of the limit
$$u_{K}(x,y)=\lim_{n \to \infty} \frac1n \sum_{k = 0}^{n-1} {\bf 1}_{B(0, K)}(\Phi_k(x) +
y).$$ The function $u_{K}$ is integrable on $X \times \R^d$,
nonnegative and $T_\Phi$-invariant. Suppose that $u_{K} \ne 0$ on a
set of positive measure. Then $u_{K} \lambda$ is a finite $T_\Phi
$-invariant measure on $X \times \R^d$, absolutely continuous with
respect to $\lambda$. Since $(X, T, \mu)$ is ergodic, this implies
that $\Phi$ is a coboundary \cite{Co79}, contrary to the assumption.
Therefore $u_{K}= 0$ a.e. for the measure $\lambda$.

Taking $K= M+1$, since ${\bf 1}_{B(0, M+1)}(\Phi_k(x) + y) \geq {\bf
1}_{B(0, M)}(\Phi_k(x))$, for $\|y\| \leq 1$, this implies:
$$\lim_{n \to \infty} \frac1n \sum_{k = 0}^{n-1} 1_{B(0,
M)}(\Phi_k(x))=0, \, \text{for a.e. } x.$$

Now we compare the discrete sum with the integral:
$$\int_0^1 {\bf 1}_{B(0, M)}(W_n(x,s)) \, ds= \frac1n \sum_{k=0}^{n-1} \int_k^{k+1}
{\bf 1}_{B(0, M)} (\Phi_k(x) + (t-k) \Phi(T^k x)) \, dt.$$

Let $\varepsilon > 0$ be arbitrary and $K$ be such that $\mu(|\Phi|
> K) \leq \varepsilon$. We have for $t \in [k, k+1]$:
$${\bf 1}_{B(0, M)} (\Phi_k(x) + (t-k) \Phi(T^k x)) \leq
{\bf 1}_{|\Phi(T^k x)| > K} + {\bf 1}_{B(0, M+K)} (\Phi_k(x)),$$ so
that for a.e. $x$.
\begin{eqnarray*}
&&\limsup_n \int_0^1 {\bf 1}_{B(0, M)}(W_n(x,s)) \, ds \\
&&\leq \limsup_n \frac1n  \sum_{k=0}^{n-1} {\bf 1}_{B(0,M +K)}
(\Phi_k(x)) + \lim_n \frac1n  \sum_{k=0}^{n-1} {\bf 1}_{(|\Phi| >
K)}(T^k x) \leq 0 + \varepsilon.
\end{eqnarray*}
It implies (\ref{asympfreq}). \eop

\begin{nota} {\rm Let us take $a > 0$ and $\underline u$ in the
unit sphere in $\R^d$. For $a$ and $\underline u$ fixed, for every
$t > 0$ we denote by $C_{a, t}$ or simply $C_t$ the cone
$\{\underline v \in \R^d: |{ \underline v \over \| \underline v} \|
- \underline u \| < a t\}$.}
\end{nota}

\vskip 3mm
\begin{lem} \label{discont} Suppose $\Phi$ is not a coboundary.
The set of discontinuity points of the increasing function $t\mapsto
\limsup_{n\rightarrow \infty} \tau_{n,C_t}(x)$ is a.e. constant
(with respect to $x$). If $t$ does not belong to this set of
discontinuity points, then $\limsup_{n\rightarrow \infty}
\tau_{n,C_t}(x)$ is almost surely constant in $x$.
\end{lem} \proof \ Let us compare $\tau_{n,C_t}(x)$ and
$\tau_{n,C_t}(Tx)$. Take $\varepsilon >0$. For every $k \geq 1$, we
have $\Phi_k(Tx) = \Phi_k(x) - \Phi(x)+\Phi(T^kx)$ and
$\|\Phi_k(Tx)-\Phi_k(x)\|\leq 2\|\Phi\|_\infty$. There is $M$ such
that, if $\Phi_k(x)>M$ and $\Phi_k(x)\in C_t$, then $\Phi_k(Tx)\in
C_{t + \varepsilon}$.

Therefore, for $s\in [0, 1]$, we have $W_n(Tx,s) \in
C_{t+\varepsilon}$, when $W_n(x,s) \in C_{t}$ and $W_n(Tx,s) \geq
M$. This implies:
\begin{eqnarray*}
\int_0^1{\bf 1}_{C_t}(W_n(x,s)) \, ds &=&\int_0^1{\bf 1}_{C_t\cap
B(0,M)}(W_n(x,s))  \, ds + \int_0^1{\bf 1}_{C_t\cap B(0,M)^c} (W_n(x,s))  \, ds\\
&\leq &\int_0^1{\bf 1}_{B(0,M)} (W_n(x,s)) \, ds + \int_0^1{\bf
1}_{C_{t + \varepsilon}} (W_n(Tx,s))  \, ds.
\end{eqnarray*}
When $n$ tends to infinity, the first integral tends to 0 almost
surely by (\ref{asympfreq}) (Lemma \ref{retcompact}) if $\Phi$ is
not a coboundary. It follows:
$$\limsup_{n} \tau_{n,C_t}(x)\leq \limsup_{n} \tau_{n, C_{t + \varepsilon}}(Tx).$$
In the same way, we have $\limsup_{n} \tau_{n,C_t}(Tx)\leq
\limsup_{n} \tau_{n,C_{t + \varepsilon}}(x)$. It follows, for every
positive real numbers $s<t<u<v$:
$$\limsup_{n\rightarrow \infty} \tau_{n,C_s}(x)\leq
\limsup_{n\rightarrow \infty} \tau_{n,C_t}(Tx)\leq
\limsup_{n\rightarrow \infty} \tau_{n,C_u}(x)\leq
\limsup_{n\rightarrow \infty} \tau_{n,C_v}(Tx).$$

This implies $\lim_{s\rightarrow t, s>t} \limsup_{n\rightarrow
\infty} \tau_{n,C_s}(x) = \lim_{s\rightarrow t,s>t}
\limsup_{n\rightarrow \infty} \tau_{n,C_s}(Tx)$, for every $t$.

Thus the common limit defines a map from $X$ into the set of
increasing functions on $[0,1]$. This map is $T$-invariant, hence
almost surely constant by ergodicity of $T$. In particular the
finite or countable set of discontinuity points of $t\mapsto
\limsup_{n\rightarrow \infty} \tau_{n,C_t}(x)$ is independent of
$x$. Outside this at most countable set of values of $t$,
$\limsup_{n\rightarrow \infty} \tau_{n,C_t}(x)$ does not depend on
$x$. \eop

Let $A_0(C) := \{\chi \in \Cal C_d([0,1]): \int_0^1 {\bf
1}_{\partial C} (\chi(s)) \ ds > 0 \}$ be the set of functions
taking their values in the boundary of $C$ for a set of positive
measure of the variable $s$.

\begin{lem} \label{Wien0} If the Lebesgue measure of the boundary
of $C$ is zero, then the Wiener measure of the set $A_0(C) $ is 0.
\end{lem} \proof \ \ The Wiener measure of $A_0(C)$ is
\begin{eqnarray*}
W(A_0(C)) = \PP((B_.)\in A_0(C)) = \PP(\{\omega : \int_0^1 {\bf
1}_{\partial C} (B_s(\omega)) \ ds > 0 \}).
\end{eqnarray*}
From the assumption on $C$, we have $\PP(B_s \in \partial C)=0$ for
every $s$ and therefore $\displaystyle {\mathbb{E}\left(\int_0^1
{\bf 1}_{\partial C}(B_s) \ ds \right) =\int_0^1\PP(B_s \in \partial
C) \, ds = 0}$. \eop

\begin{rem} \label{remAtom} It is clear that the set $\Delta$ of atoms of the
distribution of $\int_0^1{\bf 1}_C(B_s) \, ds$ (the image
probability on $[0,1]$ of the Wiener measure on $\Cal C_d([0,1])$)
is at most countable. If $c \not  \in \Delta$, then the set of
functions $\chi$ in $\Cal C_d([0,1])$ such that $\int_0^1{\bf
1}_C(\chi(s)) \, ds = c$ has zero measure for the Wiener measure.
\end{rem}

For $\eta>0$, denote by $\partial C(\eta)$ the set of points in
${\mathbb R}^d$ at a distance $\leq \eta$ from the boundary of $C$.

\vskip 3mm {\bf Proof of Theorem \ref{invPrincThm}}: Since the
interior of $C$ is non empty, there is a family of cones $C_{t} =
C_{a,t}$ contained in $C$, constructed like the cones introduced in
the preceding lemmas. We take a real $t>0$ such that
$\limsup_{n\rightarrow \infty} \tau_{n,C_t(x)}$ is almost surely
constant in $x$. If we show that $ \limsup_{n\rightarrow \infty}
\tau_{n,C_t}=1$, then also $\limsup_{n\rightarrow \infty}
\tau_{n,C}=1$.

From now, on we replace $C$ by $C_t$ still denoted $C$. In
particular, the boundary of $C$ now has Lebesgue measure 0. The
invariance principle "$W_n(.)/\sqrt{n}\rightarrow B_.$" means that,
for every continuous functional $F$ on $\Cal C_d([0,1])$, we have
\begin{eqnarray}
\mathbb{E}(F({W_n(*, .) \over \sqrt{n}}))\rightarrow
\mathbb{E}(F(B_.)).\label{invPrincForm}
\end{eqnarray}

Suppose that a sequence of probability measures $({\mathbb P}_n)$
defined on a space $X$ converges weakly to ${\mathbb P}$. By Theorem
2.7 in Billingsley's book \cite{Bi}, if a measurable function $\Psi$
from $X$ to a metric space $Y$ has a set of discontinuity points of
measure zero for ${\mathbb P}$, then the sequence of pushforward
measures $({\mathbb P_{n,\Psi}})$ converges to the pushforward
measure ${\mathbb P_\psi}$.

 For $\Psi$ we take here the function $F_c$, $c > 0$, defined on the metric
space of continuous functions from $[0,1]$ to ${\mathbb R}^d$ with
the uniform norm by
\begin{eqnarray}
F_c(\chi) = {\bf 1}_{[c,+\infty[}(\int_0^1{\bf 1}_C(\chi(s)) \, ds).
\label{phic}
\end{eqnarray}

In order to apply to $F_c$ the quoted theorem and the convergence
(\ref{invPrincForm}), we have to show that the set of discontinuity
points of $F_c$ has measure zero for the Wiener measure.

Assume that $c \not \in \Delta$ (i.e., $c$ is not an atom of the
distribution of $\int_0^1{\bf 1}_C(B_s) \, ds$). Let us consider the
set ${\Cal G}_c$ of functions $\chi_0$ such that $\chi_0 \not \in
A_0(C)$ (i.e. $\{t: \chi_0(t) \in \partial C\}$ has Lebesgue measure
0) and $\int_0^1{\bf 1}_C(\chi_0(s)) \, ds \not = c$. It has full
Wiener measure by Lemma \ref{Wien0} and Remark \ref{remAtom}.

Let us show that $F_c$ is continuous on the set ${\Cal G}_c$.

Let $\varepsilon$ be such that $0 < \varepsilon <|\int_0^1 {\bf
1}_C(\chi_0(s)) \, ds-c|$. The measure of the set of times $s$ for
which $\chi_0(s)$ is at a distance less than $\eta$ from the
boundary of $C$ tends to 0 when $\eta$ tends to 0. We can take
$\eta>0$ such that this measure is less then $\varepsilon$.

Let $\chi$ be another function at a uniform distance less than
$\eta$ from $\chi_0$. If $\chi_0(s)$ is not in $\partial C(\eta)$,
then $\chi_0(s)$ and $\chi(s)$ are either both in $C^c$ or both in
$C$. Thus, we have
\begin{eqnarray*}
&&|\int_0^1{\bf 1}_C(\chi_0(s)) \, ds - \int_0^1{\bf 1}_C(\chi(s))
\, ds| \\ &&\leq  [\int_0^1{\bf 1}_{\partial C(\eta)}(\chi_0(s) +
{\bf 1}_{(\partial C(\eta))^c}(\chi_0(s))] \, |{\bf
1}_C(\chi_0(s))-{\bf 1}_C(\chi(s))| \, ds \\
&&\leq \int_0^1{\bf 1}_{\partial C(\eta)}(\chi_0(s)) \, ds +
\int_0^1{\bf 1}_{(\partial C(\eta))^c}(\chi_0(s)) \, |{\bf
1}_C(\chi_0(s))-{\bf 1}_C(\chi(s))| \, ds \leq \varepsilon + 0.
\end{eqnarray*}

Therefore: $\displaystyle {\bf 1}_{[c,+\infty[}(\int_0^1{\bf
1}_C(\chi(s)) \, ds)= {\bf 1}_{[c,+\infty[}(\int_0^1{\bf
1}_C(\chi_0(s)) \, ds)$ and we have proved that the functional $F_c$
is continuous at $\chi_0$.

 Finally we have shown, for every $c \not \in \Delta$, the
continuity of $F_c$ on the set ${\Cal G}_c$ which has full Wiener
measure.

 For $c$ outside $\Delta$ (which is at most countable), it follows
from the theorem mentioned above:
$$\E({F_c(W_n(*,.)) \over \sqrt{n}}) \rightarrow
\mathbb{E}(F_c(B_.)),$$ that is:
\begin{eqnarray}
\lim_n \PP([\int_0^1{\bf 1}_C({W_n(x,s) \over \sqrt{n}}) \, ds ]
\geq c) = \PP([\int_0^1{\bf 1}_C(B_{s}) \, ds] \geq c).
\label{limit}
\end{eqnarray}
As $C$ is a cone, we have:
$$\int_0^1{\bf 1}_C({W_n(x,s) \over \sqrt{n}}) \, ds= \tau_{n,C}(x).$$

Because of (\ref{Palpha}), (\ref{limit}) implies that, for every
$c<1$ with $c \not \in \Delta$, $\PP(\tau_{n,C}\geq c)>0$ for $n$
large enough. As a consequence, we have
$$\PP(\limsup_n \tau_{n,C}\geq c)\geq \limsup_n\PP(\tau_{n,C} \geq c)>0.$$
Hence, $\limsup_n \tau_{n,C}$ being constant, it follows $\limsup_n
\tau_{n,C}\geq 1$. As $\tau_{n,C} \in [0,1]$, this proves $\limsup_n
\tau_{n,C}=1$. \eop

\begin{rem} We can also consider the piecewise
constant function: $V_{n}(x,s):= \Phi_k(x)\ {\rm for} \ s\in [k/n,
(k+1)/n[$. If $\Phi$ is bounded, then
$$\|W_n(x,\cdot)-V_n(x,\cdot)\|_\infty\leq \|\Phi\|_\infty.$$
On the other hand, we have
$$\int_0^1{\bf 1}_C(V_n(x,s)) \, ds={{1}\over{n}}\int_0^n{\bf 1}_C(V_n(x,{{t}\over{n}}))dt
=  \, \frac1n \Card\{k \leq n: \ \Phi_k(x)\in C\}.$$

Reasoning as before we can show that if a cone $C$ contains a cone
of the form $C_{t}$, like in Lemma \ref{discont}, then
$$\limsup_n \int_0^1{\bf 1}_C(V_n(x,s)) \, ds\geq
\limsup_n\int_0^1{\bf 1}_{C_{t}}(W_n(x,s)) \, ds=1.$$ This means
that, if $\Phi$ is a bounded function satisfying Donsker's
invariance principle, we also have the following discrete version of
the property claimed in the theorem:
$$\limsup_n \, \frac1n \Card\{k \leq n: \ \Phi_k(x)\in C\}=1.$$
\end{rem}

\vskip 3mm {\bf Acknowledgement} The authors thank B. Delyon and S.
Gouëzel for fruitful discussions. This paper is submitted to the
proceedings of the workshop in Ergodic Theory and Dynamical Systems
held at the Department of Mathematics UNC Chapel Hill in 2013-2014.
The authors are grateful to the referee for his careful reading.

\end{document}